\documentclass[final,3p,11pt, times]{article}
\usepackage{authblk}
\usepackage{amsmath}
\usepackage{amssymb,latexsym}
\usepackage[german,english]{babel}
\usepackage{url}
\usepackage{graphicx}
\usepackage{gastex}
\usepackage{longtable}
\usepackage{lscape}
\usepackage{tabularx}
\usepackage{multicol}
\usepackage{verbatim}
\usepackage{multirow}
\usepackage{graphicx}
\usepackage{epsfig,graphics,graphicx}
\usepackage{geometry}
\usepackage[usenames]{xcolor}
\usepackage{authblk}
\usepackage{verbatim}
\usepackage{multirow}
\usepackage{ragged2e}
\usepackage{xcolor}
\usepackage[T1]{fontenc}
\usepackage{authblk}
\usepackage{lipsum}
\usepackage{hyperref}
\usepackage{tikz}

\usepackage{lineno}

\hbadness=\maxdimen
\newtheorem{eg}{{\bf Example}}[subsection]

\newtheorem{deff}{Definition}[subsection]
\newtheorem{lem}{{\bf Lemma}}[subsection]
\newtheorem{thm}{{\bf Theorem}}[section]

\newtheorem{ques}{Question}

\begin{document}
	
	\title{{On the power graph of a certain gyrogroup}}

	\author[1]{Yogendra Singh}
	\author[1]{Anand Kumar Tiwari}
	\author[2]{Fawad Ali}
	\author[3]{Mani Shankar Pandey}


	\affil[1]{\small Department of Applied Science, Indian Institute of Information Technology, Allahabad 211015, India}
	
	\affil[2]{\small Institute of Numerical Sciences, Kohat University of Science and Technology, Kohat 26000, Pakistan}

	\affil[3]{\small The Institute of Mathematical Sciences, Chennai 600113, India \hspace{8cm} E-mail: yogendraiiita@gmail.com, anand@iiita.ac.in, fawadali@math.qau.edu.pk, manishankarpandey4@gamil.com}

	\maketitle
	
	\hrule
	
	\begin{abstract} 
The power graph $P(G)$ of a group $G$ is a simple graph
with the vertex set $G$ such that two distinct vertices $u,v \in G$ are adjacent in $P(G)$ if and only if $u^m = v$ or $v^m = u$, for some $m \in \mathbb{N}$. The purpose of this paper is to introduce the notion of a power graph for gyrogroups. Using this, we investigate the combinatorial properties of a certain gyrogroup, say $G(n)$, of order $2^n$ for $n \geq 3$. In particular, we determine the Hamiltonicity and planarity of the power graph of $G(n)$. Consequently, we calculate distant properties, resolving polynomial, Hosoya and reciprocal Hosoya polynomials, characteristic polynomials, and the spectral radius of the power graph of $G(n)$.

\end{abstract}

\smallskip



\hrule

\section{Introduction}

Gyrogroup, a non-associative algebraic structure that is a generalisation of a group, was first described by Ungar in 1988, see  \cite{AU}. The concept of gyrogroup initially appeared in the study of Einstein's relativistic velocity addition law, see \cite{ungar_thomas,ungar_book}. In this direction, Sukumran with Wiboonton \cite{suksumran_gyroaction,TS, suksumran_lagrange,suksumran_isomorphism} studied several well-known group theoretic results for gyrogroups such as the Lagrange theorem, the fundamental isomorphism theorem, and the Cayley theorem, etc. Alike groups, Foguel and Ungar \cite{fu(2001)} studied the concept of transversal in gyrogroups and established a relationship between twisted subgroups of a group and gyrogroups. Bussaban et al. \cite{Bussaban(2019)} studied the combinatorial properties of gyrogroups for the first time by introducing the Cayley graphs for gyrogroups. For more recent developments and in-depth knowledge of the Cayley graph associated with gyrogroups, see \cite{Maungchang(2021),Maungchang(2022)}. The notion of the power graph of a group is a relatively new addition to the domain of graphs associated with groups. Chakrabarty et al. \cite{Chakrabarty(2009)} described the concept of undirected {\it power graph} for a group for the first time. Recently, Ajay et al. presented almost all the existing results on power graphs of a group in the survey article \cite{ALPT}.

Since gyrogroups are a generalisation of groups, it is natural to explore its combinatorial properties. In order to do this, we define the power graph of a gyrogroup (as defined for a group) as follows:

\begin{deff}
	The power graph $P(G)$ of a gyrogroup $G$ is the simple undirected graph with vertex set $G$ and two distinct vertices $u,v \in G$ are adjacent in $P(G)$ if and only if $u^m = v$ or $v^m = u$, for some $m \in \mathbb{N}$.	
\end{deff}

In this research, we investigate combinatorial structures of a gyrogroup $G(n)$ of order $2^n$, $n \geq 3$. In particular, we study Hamiltonicity and planarity of $G(n)$. In addition, we calculate Hosoya and the reciprocal Hosoya polynomials, resolving polynomial, characteristic polynomials and spectral radius of $G(n)$. Also, we calculate distant as well as detour distant properties such as closure, interior, distance degree sequence of power graphs of $G(n)$.


\section{Preliminaries}
In this section, we present basic knowledge of graphs and gyrogroups. Most of the following graph theory related terms are given in \cite{Ali(2020)}. Also, we refer to \cite{West(2009)}, \cite{AL(2019)} for their detailed description.  

Throughout this article, by a graph $G$ we mean a simple graph, i.e., an undirected graph without loops or multiple edges. The sets $V(G)$ and $E(G)$ represent the vertex set and edge set of $G$, respectively. The order $|G|$ of a graph $G$ is the cardinality of $V(G)$. Let $G_1$ and $G_2$ be two graphs with the vertex sets $V(G_1)$ and $V(G_2)$ respectively. Then the union $G_1 \cup G_2$ of $G_1$ and $G_2$ is a graph with the vertex set $V(G_1 \cup G_2) = V(G_1) \cup V(G_2)$ and the edge set $E(G_1 \cup G_2) = E(G_1) \cup V(G_2)$. The graphs $G_1$ and $G_2$ are isomorphic if there is a bijective map between their vertex sets $V(G_1)$ and $V(G_2)$ that sends vertices to vertices, edges to edges, and  preserves incidence.

Let $S \subseteq V(G)$. Then the induced subgraph $G[S]$ of $G$ is the graph with the vertex set $S$ and the edge set consists of all the edges in $E(G)$ which have both of its endpoints in $S.$ If two distinct vertices $u, v$ are adjacent by an edge, then we write $u \sim v$; otherwise, $u \nsim v$. In a non-trivial connected graph $G$, the distance $d(u,v)$ between two vertices $u$ and $v$ is the length of a shortest $u$-$v$ path in $G$. While, the length of a largest $u$-$v$ path between $u$ and $v$ is known as detour distance and is denoted by $d_D(u,v)$. The eccentricity (detour eccentricity) of a vertex in $G$ is the largest  distance (detour distance) between a vertex $u$ and any other vertex of $G$, and it is denoted by $e(u) \, (e_D(u))$. The largest eccentricity (detour eccentricity) among all the vertices in $G$ is the diameter (detour diameter) of $G$ and is denoted by  $dia(G)\, (dia_D(G)) $. Also, the smallest eccentricity (detour eccentricity) among all the vertices in $G$ is the radius (detour radius) of $G$ and is denoted by  $rad(G) \,(rad_D(G)).$ 


For a vertex $u \in V(G)$, the set $N(u) = \{ v \in G \ | \ u \sim v\}$ is known as the open neighborhood, while $N[u] = N(u) \cup \{v\}$ is the closed neighborhood of $u$. The degree $\deg(u)$ of the vertex $u$ is defined as the cardinality of the set $N(u)$. If $N[u] = N[v]$ or $N(u) = N(v)$ holds for every $u, v \in V(G)$, $u \neq v$, then $u$ and $v$ are said to be adjacent or non-adjacent twins, respectively. In general, we say that $u$ and $v$ are twins if they are adjacent or non-adjacent twins. A subset $U \subseteq V(G)$ is called a twin-set in $G$, if every distinct pair of vertices in $U$ are twins. 

A vertex $u$ of $G$ is a boundary vertex of a vertex $v$ if $d(w,v) \leq d(u,v)$ for all $w \in N(u)$. A vertex $u$ of $G$ is a boundary vertex of $G$ if $u$ is a boundary vertex of some vertex of $G$. The boundary $ \partial(G)$ of $G$ is the subgraph of $G$ induced by its boundary vertices.  A vertex $v$ is called a center vertex of $G$ if $e(v) = rad (G)$. The center $C(G)$ of $G$ is the collection of all its center vertices. If for every vertex $u$ other than $v$, there is a vertex $w$ such that $v$ lies between $u$ and $w$, then $v$ is called an interior vertex of $G$. The subgraph of $G$ induced by interior vertices is known as the interior $Int(G)$ of $G$. 

Let $G$ be a graph of order $n$. Then the graph obtained from $G$ by adding edges between non-adjacent vertices whose degree sum is at least $n$ until this can no longer be done is known as the closure of $G$ and is denoted by $Cl(G)$. A vertex $u$ is called a complete vertex, if the induced subgraph of $N(u)$ is a complete graph.

Let $U = \{u_1, u_2, \ldots , u_s\}$ be an ordered subset of $V(G)$ and $v \in V(G)$. Then the $s$-vector $r(v|U) = (dis(v,u_1),dis(v,u_1),\ldots, dis(v,u_s))$ is called the representation of $v$ with respect to $U.$ If any two different vertices of $G$ have different representations with respect to $U$, then $U$ is said to be a resolving set for $G$. A minimal resolving set of $G$ is called a basis of $G$, and the cardinality of a basis is called the metric dimension $\psi(G)$ of $G$. 

Let $|G|=n$ and let $R(G, k)$ denote the collection of resolving sets which are $k$-subsets. If $|R(G, k)| = r_k$, then a resolving polynomial
of $G$ is defined as $\psi(G, x)= \sum_{k=\psi(G)}^{n} r_k \psi(G)$. The sequence $(r_{\psi(G)}, r_{\psi(G) +1}, \ldots , r_n)$ of coefficients of $\psi(G)$ is known as resolving sequence.

For $u \in  V(G)$, let $deg_k(u)(Deg_k(u))$ be the set of vertices at distance (detour distance) $k$ from $u$. Then
$(deg_0(u)$, $deg_1(u)$, $deg_2(u),\ldots, deg_{ec(u)}(u)((Deg_0(u)$, $Deg_1(u)$, $Deg_2(u), \ldots, Deg_{ecD(u)}(u))$ is known as distance (detour distance) degree sequence and is denoted by $dds(u) (dds_D(u))$. The $dds(G)\, (dds_D(G))$ of a graph $G$ is the set $dds(u)\, (dds_D(u))$, where $u \in V(G)$. 

A graph is called planar if it can be embedded on the plane. A cycle in a graph $G$ is called Hamiltonian if it covers all the vertices. A graph $G$ is called Hamiltonian graph if it contains a Hamiltonian cycle.

Let $G$ be a graph with the vertex set $V(G) = \{u_1, u_2, \ldots, u_n\}$. Then the $n \times n$ symmetric matrix $A_G= [a_{ij}]$ is called the adjacency matrix of $G$, where $a_{ij}=1$ if $u_i$ is adjacent with $u_j$ and 0 otherwise. The characteristic polynomial of $G$ is $p(\lambda) = \begin{vmatrix} \lambda I-A_G \end{vmatrix}$ and then the eigenvalues of $G$ are the eigenvalues of matrix $A_G$. The spectrum of $G$ is defined as the set of all eigenvalues along with its multiplicity. Since $A_G$ is symmetric, all of its eigenvalues are real. Therefore, these eigenvalues can be arranged in decreasing order as $\lambda_1(G) \geq \lambda_2(G) \geq \ldots \geq \lambda_n(G)$. The greatest eigenvalue $\lambda_1(G)$ of $G$ is known as the spectral radius.

Now, we recall the definition of a gyrogroup and some of its properties that we use further.

\begin{deff} \label{gyrogroup}  
	A pair $(G,\circ)$ consisting of a non-empty set $G$ and a binary operation $\circ$ on $G$ is called a {\it gyrogroup} if its binary operation $\circ$ satisfies the following axioms. 
	\begin{enumerate}
		\item There is at least one element $e\in G$, called a left identity, such that $$e\circ a=a\,\,\text{ for all } a\in G.$$
		
	\item For every $a\in G,$ there is an element $a^{-1}\in G$, called a left inverse of $a,$ such that $a^{-1} \circ \ a=e$.
		
	\item There exists an automorphism $gyr[a,b] \in Aut(G, \circ)$ such that 
	$$ a\circ(b\circ c)= (a\circ b)\circ gyr[a,b] c \  \text{ for all } a,b,c\in G.$$
		
   \item $gyr[a\circ b, b]=gyr[a,b]$ for all $a,b\in G.$
	\end{enumerate} 
\end{deff}

\begin{deff} A gyrogroup $(G, \circ)$ is gyro-commutative if its binary operation obeys the gyro-commutative law 
	$$ a \circ b= gyr[a,b](b \circ a) \text{ for all } a,b\in G.$$
	
\end{deff}

\begin{deff}
	A gyrogroup that is not a group is called a non-degenerate gyrogroup.
\end{deff}

\section{Power graph of a gyrogroup} \label{s3}

In \cite{MASU}, Madhavi et al. proposed the gyrogroup $G(n)$ of order $2^n$ given in Example \ref{eg1}. We present its power graph construction in Example \ref{eg2}. In this paper, we determine the combinatorial properties of the gyrogroup $G(n)$ through its power graph $P(G(n)).$

\begin{eg}\label{eg1}
	
	For a natural number $n \geq 3$, let $P(n) = \{0, 1, 2, \ldots, 2^{n-1} -1\}$, $H(n) = \{2^{n-1}, 2^{n-1} + 1, \ldots, 2^{n} - 1 \}$, and $G(n) = P(n) \cup H(n)$. Clearly, $P(n)$ is a cyclic group under addition madulo $m = 2^{n-1}.$ The binary operation of the gyrogroup $(G(n),\oplus)$ is
	defined as follows:
	
	\begin{equation*}
	i \oplus j=
	\begin{cases}
	t & \text{$(i,j) \in P(n) \times P(n)$}\\
	t + m & \text{$(i,j) \in P(n) \times H(n)$}\\
	s + m & \text{$(i,j) \in H(n) \times P(n)$} \\
	k & \text{$(i,j) \in H(n) \times H(n)$}
	\end{cases}       
	\end{equation*}
	
	where $t,s,k \in P(n)$ are the following non-negative integers:
	
	$$ \begin{cases}
	t = i + j & \text{mod m} \\
	s = i + (\frac{m}{2}-1)j & \text{mod m}\\
	k = (\frac{m}{2}+1)i + (\frac{m}{2}-1)j & \text{mod m}\\
	\end{cases} $$
	
\end{eg} 

\begin{eg}\label{eg2} Since $i^2 = i \oplus i = k = ((\frac{m}{2}+1)i + (\frac{m}{2}-1)i)) \ mod \ m = 0$ for all $(i,i) \in H(n) \times H(n)$, every element of $H(n)$ in $P(G(n))$ is adjacent with identity only. As $P(n)$ is a cyclic group of order $2^{n-1}$, by \ref{l1} its power graph is isomorphic to complete graph $K_{2^{n-1}}$. Thus the power graph of $(G(n),\oplus)$ is $K_{2^{n-1}} \cup \{0-2^{n-1}\} \cup \{0-{2^{n-1}+1}\} \cup \cdots \cup \{0-{2^{n-1}+1}\}$.
	
\end{eg}

First, we discus the Hamiltonicity and planarity of $P(G(n)).$ For this, we need the following result:

\begin{thm}\label{l1} $[\cite{Chakrabarty(2009)}, Theorem \ 2.12]$ The power graph $P(G)$ of a finite group $G$ is complete if and only if $G$ is																										 a cyclic group of order $p^m$ for some prime $p$ and non-negative integer $m$.
\end{thm}

\begin{lem}
	The power graph $P(G(n))$ is planar for $n = 3$ and non-planar for $n \geq 4$.
\end{lem}
\noindent{\bf Proof.} For $n=3$, the power graph is isomorphic to the graph shown in Figure 1, which is clearly planar. Let $n \geq 4$ and $m = 2^{n-1}.$ As $P(n)$ is a cyclic group under addition madulo $m = 2^{n-1}$, by Theorem \ref{l1}, power graph of $P(n)$ is $K_{2^{n-1}}$ and hence contains a copy of $K_5$. Thus, $P(G(n))$ also contains a copy of $K_5$. Hence, $P(G(n))$ is not planer. \hfill $\Box$

\begin{lem}
	The power graph $P(G(n))$ is not Hamiltonian.
\end{lem}

\noindent{\bf Proof.} Note that $P(G(n))$ has $2^{n-1}$ vertices of degree 1. These $2^{n-1}$ vertices along with $e$ induces a tree. Hence, $P(G(n))$ is not Hamiltonian. \hfill $\Box$

\smallskip

Cameron and Ghosh \cite{CG} proved that non-isomorphic finite groups may have isomorphic power graphs, however, finite abelian groups with isomorphic power graphs are always isomorphic. By the argument given in Example \ref{e3.1}, we see that non-isomorphic finite gyrogroups may have isomorphic power graphs.

\begin{eg}\label{e3.1}
	Consider the gyrogroups $K(1)$ and $N(1)$, given in \cite{MAS}. The Cayley tables and associated gyration tables of the gyrogroups $K(1)$ and $N(1)$ are given in Tables 1 and 2, respectively.  The power graph of $K(1)$ and $N(1)$ are shown in Figures 1 and 2, respectively. Define a map $f: P(K(1)) \rightarrow P(N(1))$ such that $f(0)=0, f(1) = 1, f(2)=7, f(3)=6, f(4) = 2, f(5)=3, f(6) = 5$, and $f(7)=4$. Clearly, $f$ is an isomorphism. Thus, we have gyrogroups $K(1)$ and $N(1)$ such that  $P(K(1)) \cong P(N(1))$ but $K(1) \ncong N(1)$. 
\end{eg}	

\begin{picture}(0,0)(-50,7)
\setlength{\unitlength}{7.5mm}

\drawline[AHnb=0](0,0)(0,1)
\drawline[AHnb=0](0,0)(0,-1)
\drawline[AHnb=0](0,0)(-1,-1)
\drawline[AHnb=0](0,0)(-1,1)
\drawline[AHnb=0](0,0)(-.65,0)
\drawline[AHnb=0](-.65,0)(-1,1)
\drawline[AHnb=0](-.65,0)(-1,-1)
\drawline[AHnb=0](-1,1)(-1,-1)
\drawline[AHnb=0](0,0)(1,1)
\drawline[AHnb=0](0,0)(1,-1)
\put(.15,-.1){\scriptsize {\tiny $0$}}
\put(.15,1){\scriptsize {\tiny $2$}}
\put(.15,-1){\scriptsize {\tiny $5$}}
\put(-.65,.03){\scriptsize {\tiny $6$}}
\put(-1.25,-1){\scriptsize {\tiny $7$}}
\put(-1.25,1){\scriptsize {\tiny $1$}}
\put(1.1,1){\scriptsize {\tiny $3$}}
\put(1.1,-1){\scriptsize {\tiny $4$}}

\put(-2,-2){\scriptsize {\tiny {\bf Figure 1:} Power graph of $K$}}
\end{picture}

\begin{picture}(0,0)(-105,2)
\setlength{\unitlength}{7.5mm}

\drawline[AHnb=0](0,0)(0,1)
\drawline[AHnb=0](0,0)(0,-1)
\drawline[AHnb=0](0,0)(-1,-1)
\drawline[AHnb=0](0,0)(-1,1)
\drawline[AHnb=0](0,0)(-.75,0)
\drawline[AHnb=0](-.75,0)(-1,1)
\drawline[AHnb=0](-.75,0)(-1,-1)
\drawline[AHnb=0](-1,1)(-1,-1)
\drawline[AHnb=0](0,0)(1,1)
\drawline[AHnb=0](0,0)(1,-1)

\put(.15,-.1){\scriptsize {\tiny $0$}}
\put(.15,1){\scriptsize {\tiny $7$}}
\put(.15,-1){\scriptsize {\tiny $3$}}
\put(-.65,.03){\scriptsize {\tiny $5$}}
\put(-1.25,-1){\scriptsize {\tiny $4$}}
\put(-1.25,1){\scriptsize {\tiny $1$}}
\put(1.1,1){\scriptsize {\tiny $6$}}
\put(1.1,-1){\scriptsize {\tiny $2$}}

\put(-2.5,-2){\scriptsize {\tiny {\bf Figure 2:} Power graph of $N$}}

\end{picture}

\newpage
\hspace{7cm}	\textbf{Table 1} 
\begin{center}
	
	\begin{tabular}{c | c c c c c c c c ||c | c c c c c c c c}
		
		$\oplus_K$ & 0 & 1 & 2 & 3 & 4 & 5 & 6 & 7 & $gyr_K$ & 0 & 1 & 2 & 3 & 4 & 5 & 6 & 7\\
		\cline{1-18}

		0 & 0 & 1 & 2 & 3 & 4 & 5 & 6 & 7 & 0 & I & I & I & I & I & I & I & I \\
		1 & 1 & 0 & 3 & 2 & 5 & 4 & 7 & 6 & 1 & I & I & I & I & I & I & I & I \\
		2 & 2 & 3 & 0 & 1 & 6 & 7 & 4 & 5 & 2 & I & I & I & I & A & A & A & A \\
		3 & 3 & 2 & 1 & 0 & 7 & 6 & 5 & 4 & 3 & I & I & I  & I  &A & A& A& A \\

		4 & 4 & 5 & 6 & 7 & 0 & 1 & 2 & 3 & 4 & I & I & A & A & I & I & A & A \\
		5 & 5 & 4 & 7 & 6 & 1 & 0 & 3 & 2 & 5 & I & I & A & A & I & I & A & A \\
		6 & 6 & 7 & 4 & 5 & 3 & 2 & 1 & 0 & 6 & I & I & A & A & A & A & I & 
		I \\
		
		7 & 7 & 6 & 5 & 4 & 2 & 3 & 0 & 1 & 7 & I & I & A & A & A & A & I & I
		
	\end{tabular}
\end{center}

\bigskip

\hspace{7cm}  \textbf{Table 2} 
\begin{center}

	\noindent\begin{tabular}{c | c c c c c c c c ||c | c c c c c c c c}
		
		$ \oplus_N$	& 0 & 1 & 2 & 3 & 4 & 5 & 6 & 7 & $gyr_N$ & 0 & 1 & 2 & 3 & 4 & 5 & 6 & 7\\
		\cline{1-18}
		
		0 & 0 & 1 & 2 & 3 & 4 & 5 & 6 & 7 & 0 & I & I  & I & I & I & I & I & I \\
		1 & 1 & 0 & 3 & 2 & 5 & 4 & 7 & 6 & 1 & I & I  &I & I & I & I & I & I \\
		2 & 2 & 3 & 0 & 1 & 6 & 7 & 4 & 5 & 2 & I & I & I & I & D & D & D & D \\
		3 & 3 & 2 & 1 & 0 & 7 & 6 & 5 & 4 & 3 & I & I  & I & I & D & D & D & D \\
		4 & 4 & 5 & 6 & 7 & 1 & 0 & 3 & 2 & 4 & I & I & D & D & I & I & D & D \\
		5 & 5 & 4 & 7 & 6 & 0 & 1 & 2 & 3 & 5 & I & I & D & D & I & I & D & D \\
		6 & 6 & 7 & 5 & 4 & 3 & 2 & 0 & 1 & 6 & I & I  & D & D & D & D & I & I \\
		7 & 7 & 6 & 4 & 5 & 2 & 3  & 1 & 0 & 7 & I & I & D & D & D & D & I & I \\
		
	\end{tabular}
\end{center}

Now, one may ask naturally whether finite gyro-commutative gyrogroups with isomorphic power graphs are isomorphic. By the Example \ref{e3.2}, we see that it need not be true.    

\begin{eg} \label{e3.2}
	Consider the gyro-commutative gyrogroups $G_8$ \cite{TS} and $M(1)$ \cite{MAS}. The Cayley tables and associated gyration tables of the gyrogroups $G_8$ and $M(1)$ are given in Tables 3 and 4, respectively. Their power graphs are shown in Figures 3 and 4, respectively. Define a map $f: P(G_8) \rightarrow P(M(1))$ such that $f(0)=0, f(1) = 3, f(2)=7, f(3)=5, f(4) = 4, f(5)=6, f(6) = 1$, and $f(7)=2$. Clearly, $f$ is an isomorphism. Thus, $P(G_8) \cong P(M(1))$ but $G_8 \ncong M(1)$. 

\bigskip

	\hspace{7cm}	\textbf{Table 3} 
	\begin{center}
		
		\begin{tabular}{c | c c c c c c c c ||c | c c c c c c c c}
			
			$\oplus$ & 0 & 1 & 2 & 3 & 4 & 5 & 6 & 7 & $gyr$ & 0 & 1 & 2 & 3 & 4 & 5 & 6 & 7\\
			\cline{1-18}
			
			0 & 0 & 1 & 2 & 3 & 4 & 5 & 6 & 7 & 0 & I & I & I & I & I & I & I & I \\
			1 & 1 & 3 & 0 & 2 & 7 & 4 & 5 & 6 & 1 & I & I & I & I & A & A & A & A \\
			2 & 2 & 0 & 3 & 1 & 5 & 6 & 7 & 4 & 2 & I & I & I & I & A & A & A & A \\
			3 & 3 & 2 & 1 & 0 & 6 & 7 & 4 & 5 & 3 & I & I & I  & I & I & I& I& I \\

			4 & 4 & 5 & 7 & 6 & 3 & 2 & 0 & 1 & 4 & I & A & A & I & I & A & I & A \\
			5 & 5 & 6 & 4 & 7 & 2 & 0 & 1 & 3 & 5 & I & A & A & I & A & I & A & I \\
			6 & 6 & 7 & 5 & 4 & 0 & 1 & 3 & 2 & 6 & I & A & A & I & I & A & I & 
			A \\
			
			7 & 7 & 4 & 6 & 5 & 1 & 3 & 2 & 0 & 7 & I & A & A & I & A & I & A & I

		\end{tabular}
		
	\end{center}

	\hspace{6.5cm}	\textbf{Table 4} 
	\begin{center}

		\noindent\begin{tabular}{c | c c c c c c c c ||c | c c c c c c c c}
			
			$ \oplus_M$	& 0 & 1 & 2 & 3 & 4 & 5 & 6 & 7 & $gyr_M$ & 0 & 1 & 2 & 3 & 4 & 5 & 6 & 7\\
			\cline{1-18}
			
			0 & 0 & 1 & 2 & 3 & 4 & 5 & 6 & 7 & 0 & I & I  & I & I & I & I & I & I \\
			1 & 1 & 0 & 3 & 2 & 5 & 4 & 7 & 6 & 1 & I & I  &I & I & I & I & I & I \\
			2 & 2 & 3 & 0 & 1 & 6 & 7 & 4 & 5 & 2 & I & I & I & I & C & C & C & C \\
			3 & 3 & 2 & 1 & 0 & 7 & 6 & 5 & 4 & 3 & I & I  & I & I & C & C & C & C \\
			4 & 4 & 5 & 6 & 7 & 1 & 0 & 3 & 2 & 4 & I & I & C & C & I & I & C & C \\
			5 & 5 & 4 & 7 & 6 & 0 & 1 & 2 & 3 & 5 & I & I & C & C & I & I & C & C \\
			6 & 6 & 7 & 5 & 4 & 2 & 3 & 1 & 0 & 6 & I & I  & C & C & C & C & I & I \\
			7 & 7 & 6 & 4 & 5 & 3 & 2  & 0 & 1 & 7 & I & I & C & C & C & C & I & I \\
			
		\end{tabular}
		
	\end{center}
	
	\newpage

	\begin{picture}(0,0)(-50,16)
	\setlength{\unitlength}{7.5mm}
	
	\drawline[AHnb=0](-1,1)(1,1)
	\drawline[AHnb=0](-1,1)(0,.75)
	\drawline[AHnb=0](0,.75)(1,1)
	\drawline[AHnb=0](0,0)(0,.75)
	\drawline[AHnb=0](0,0)(0,-1)
	\drawline[AHnb=0](0,0)(-1,-1)
	\drawline[AHnb=0](0,0)(-1,1)
	\drawline[AHnb=0](0,0)(-.65,0)
	\drawline[AHnb=0](-.65,0)(-1,1)
	\drawline[AHnb=0](-.65,0)(-1,-1)
	\drawline[AHnb=0](-1,1)(-1,-1)
	\drawline[AHnb=0](0,0)(1,1)
	\drawline[AHnb=0](0,0)(1,-1)
	\put(.15,-.1){\scriptsize {\tiny $0$}}
	\put(.15,.5){\scriptsize {\tiny $6$}}
	\put(.15,-1){\scriptsize {\tiny $3$}}
	\put(-.65,.03){\scriptsize {\tiny $4$}}
	\put(-1.25,-1){\scriptsize {\tiny $5$}}
	\put(-1.25,1){\scriptsize {\tiny $1$}}
	\put(1.1,1){\scriptsize {\tiny $7$}}
	\put(1.1,-1){\scriptsize {\tiny $2$}}

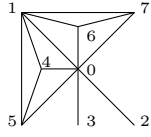
	\put(-2,-2){\scriptsize {\tiny {\bf Figure 3:} Power graph of $G_8$}}
	\end{picture}

	\begin{picture}(0,0)(-110,12)
	\setlength{\unitlength}{7.5mm}
	
	\drawline[AHnb=0](-1,1)(1,1)
	\drawline[AHnb=0](-1,1)(0,.75)
	\drawline[AHnb=0](0,.75)(1,1)
	\drawline[AHnb=0](0,0)(0,.75)
	\drawline[AHnb=0](0,0)(0,-1)
	\drawline[AHnb=0](0,0)(-1,-1)
	\drawline[AHnb=0](0,0)(-1,1)
	\drawline[AHnb=0](0,0)(-.65,0)
	\drawline[AHnb=0](-.65,0)(-1,1)
	\drawline[AHnb=0](-.65,0)(-1,-1)
	\drawline[AHnb=0](-1,1)(-1,-1)
	\drawline[AHnb=0](0,0)(1,1)
	\drawline[AHnb=0](0,0)(1,-1)
	\put(.15,-.1){\scriptsize {\tiny $0$}}
	\put(.15,.5){\scriptsize {\tiny $1$}}
	\put(.15,-1){\scriptsize {\tiny $5$}}
	\put(-.65,.03){\scriptsize {\tiny $4$}}
	\put(-1.25,-1){\scriptsize {\tiny $6$}}
	\put(-1.25,1){\scriptsize {\tiny $3$}}
	\put(1.1,1){\scriptsize {\tiny $2$}}
	\put(1.1,-1){\scriptsize {\tiny $7$}}

	\put(-2,-2){\scriptsize {\tiny {\bf Figure 4:} Power graph of $M$}}
	\end{picture}
\end{eg}

\vspace{2.35cm}

\section{Hosoya and reciprocal status Hosoya polynomial }

In 1988, Hosoya \cite{Hosoya(1971)} introduced a distance-based polynomial of a connected graph $G$, which is denoted by $H (G, x)$ and is defined as follows:

$$H(G, x) = \sum_{i \geq 0} dis(G, i)x^i.$$

The coefficient $dis(G, i)$ denotes the number of vertex pairs $(v,w)$ required for $dis(v,w) = i$, where $i \leq diam(G)$. In \cite{Ramane(2019)}, the authors introduced the reciprocal status Hosoya polynomial of $G$, which is expressed as: 

$$Hrs(G, x) = \sum_{uv \in E(G)} x^{rs(u)+rs(v)},$$

where $rs(v) = \sum_{u \in V(G), v \neq u} \frac{1}{dis(u,v)}$ denotes the reciprocal status of a vertex $v$. Also, we refer to \cite{AR(2022)} for a more detailed explanation. 

\smallskip
Here, we determine the Hosoya and reciprocal status Hosoya polynomial of $P(G(n)).$ For this, we first prove the following lemma.

\begin{lem}\label{4.1}
	Let $P(G(n))$ be the power graph of $G(n)$. Then:
\end{lem}

\begin{equation*}
dis(P(G(n),i)= 
\begin{cases}
2^n, & \text{$ for \ i = 0;$}\\
\frac{2^{n-1}(2^{n-1}+1)}{2},  & \text{$ for \ i=1;$}\\
\frac{3. 2^{n-1}(2^{n-1}-1)}{2}, & \text{$ for \ i =2.$}
\end{cases}       
\end{equation*}

\noindent{\bf Proof.} If $V$ denotes the set of all pairs of vertices of $P(G(n))$ and $C(P(G(n),i)) = \{(u,v); u, v \in V(P(G(n)) : dis(u,v) = i \}$, then $|V| = 2^n + (2^n-1) + (2^n-2) + \cdots + (2^n-(2^n-1)) =    2^{n-1}(2^n+1)$ and $ dis(P(G(n),i) = |C(P(G(n),i))|$. As $diam(P(G(n))) = 2$, we need to find $dis(P(G(n)), 0)$,
$dis(P(G(n)), 1)$, and $dis(P(G(n)),2)$. Since $d(u,u) = 0 \ \text{for all} \ u \in V(P(G(n))$, we have $$C(P(G(n),0)) = \{(u,u); u \in V(P(G(n)) : dis(u,u) = 0 \} = V(P(G(n))) \text{ and } |C(P(G(n),0))| = 2^n.$$

Now, $C(P(G(n),1)) = \{(e,v); e \neq v \ \text{and} \ v \in H(n)\} \cup \{(e,v); e \neq v \ \text{and} \ v \in P(n) \setminus \{e\}\} \cup  \{(u,v); u \neq v \ \text{and} \ u, v \in P(n) \setminus \{e\} \}.$ This gives that, 

\begin{align*}
|C(P(G(n),1))| & = dis(P(G(n)), 1)\\
& = 2^{n-1} + (2^{n-1}-1) + (2^{n-1}-1)-1 + (2^{n-1}-1)-2 + \cdots  + (2^{n-1}-1)\\
& -((2^{n-1}-1)-1) \\ & =  2^{n-1} + (2^{n-1}-1) + (2^{n-1}-1)(2^{n-1}-2) - \frac{(2^{n-1}-1)(2^{n-1}-2)}{2} \\
& = \frac{2^{n-1}(2^{n-1}+1)}{2}.
\end{align*}

Clearly, 
\begin{equation} \label{eq1}
V = C(P(G(n),0)) \cup C(P(G(n),1)) \cup C(P(G(n),2))
\end{equation}

From Equation \ref{eq1}, we have

$$|V| = dis(P(G(n)), 0) + dis(P(G(n)), 1) + dis(P(G(n)), 2)$$

Hence, 
\begin{equation*} 
\begin{split}
dis(P(G(n)), 2) & = |V| - dis(P(G(n)), 0) - dis(P(G(n)), 1)  \\
& = 2^{n-1}(2^n+1) - 2^n - \frac{2^{n-1}(2^{n-1}+1)}{2} 
\\
& = \frac{3. 2^{n-1}(2^{n-1}-1)}{2}.
\end{split}
\end{equation*}  \hfill $\Box$

\begin{lem}\label{4.2}
	Let $P(G(n))$ be the power graph of $G(n)$. Then: 
	$$H(P(G(n)), x)= 2^nx^0 + \frac{2^{n-1}(2^{n-1}+1)}{2}x^1
	+ \frac{3. 2^{n-1}(2^{n-1}-1)}{2}x^2  $$
\end{lem}

\noindent{\bf Proof.} Since 

\begin{equation}\label{eq3}
H(G, x) = \sum_{i \geq 0} dis(G, i)x^i, \ \text{where} \ i \leq diam(G)=2. 
\end{equation}

Therefore, by Lemma \ref{4.1} and Equation \ref{eq3}, we have

\begin{equation*} 
\begin{split}
H(P(G(n)), x) & = dis(P(G(n)), 0) x^0 + dis(P(G(n)), 1)x^1
+ dis(P(G(n)), 2)x^2  \\
& = 2^nx^0 + \frac{2^{n-1}(2^{n-1}+1)}{2}x^1
+ \frac{3. 2^{n-1}(2^{n-1}-1)}{2}x^2. 
\end{split} 
\end{equation*} \hfill $\Box$

\begin{lem}\label{4.3}
	Let $P(G(n))$ be the power graph of $G(n)$. Then: 
	$$H_{rs}(P(G(n)))= (2^{n-1}-1)x^{\frac{2^{n+2}-2^{n-1}-4}{2}} + 2^{n-1}x^{2^n+2^{n-1}-1} + {2^{n-1} -1 \choose 2}x^{2^n+2^{n-1}-2}
	$$
\end{lem}

\noindent{\bf Proof.} The reciprocal status Hosoya polynomial of $G$ is: 

$$H_{rs}(G, x) = \sum_{uv \in E(G)} x^{rs(u)+rs(v)}, \ \text{where} \ rs(v) = \sum_{u \in V(G), v \neq u} \frac{1}{dis(u,v)}.$$ Note that, there are three types of edges in $P(G(n))$, namely $eu$, $ev$, $uw$, where $u, w \in P(n) \setminus \{e\}$, $v \in H(n).$ Hence,

\begin{equation} \label{eq2}
H_{rs}(G, x) = \sum_{eu \in E(G)} x^{rs(e)+rs(u)} + \sum_{ev \in E(G)} x^{rs(e)+rs(v)} + \sum_{uw \in E(G)} x^{rs(u)+rs(w)}.
\end{equation}

Since $diam(P(G))=2$, we have $deg(u)$ number of vertices at distance $1$ form $u$ and $2^n-1-deg(u)$ vertices at distance 2. Hence $rs(u) =\frac{1}{1}.d(u) + \frac{1}{2}.(2^n-1-d(u))$. As degree of $u, w \in P(n) \setminus \{e\}$, $v \in H(n), \text{and} \ e$ are $2^{n-1}-1, 2^{n-1}, \text{and} \ 2^n-1$, respectively. We have  
$$rs(e)+rs(u) = \frac{1}{1}.(2^{n}-1) + \frac{1}{2}.(2^n-1-2^{n}+1) + \frac{1}{1}.(2^{n-1}-1) + \frac{1}{2}.(2^{n}-1-2^{n-1}+1)= \frac{2^{n+2}-2^{n-1}-4}{2}.$$
$$rs(e)+rs(w) = \frac{1}{1}.(2^{n}-1) + \frac{1}{2}.(2^n-1-2^{n}+1) + \frac{1}{1}.(1) + \frac{1}{2}.(2^n-1-1) =  2^n+2^{n-1}-1.$$
$$rs(u)+rs(v) = \frac{1}{1}.(2^{n-1}-1) + \frac{1}{2}.(2^{n}-1-2^{n-1}+1) + \frac{1}{1}.(2^{n-1}-1) + \frac{1}{2}.(2^{n}-1-2^{n-1}+1) =  2^n+2^{n-1}-2.$$

Now, the number of edges of $P(G(n))$ of types $eu$, $ev$, $uw$ are $2^{n-1}-1, 2^{n-1}, \text{and} {2^{n-1} -1 \choose 2}$, respectively. By using these values in Equation \ref{eq2}, we get the required proof. \hfill $\Box$

\section{Resolving polynomial of the power graph $P(G(n))$}

In this section, we find certain metric dimension properties of $P(G(n))$. Then we determine the resolving polynomial, characteristic polynomial and spectral radius of the $P(G(n))$. For this, we recall the following results:


\begin{lem} \label{l6.2} \cite{Hernando(2007)}
	Suppose the graph $G$ is connected and $u, v$ are twins in $G$, also $A$ is the resolving set of $G$. Then $u$ or $v$ is in $A$. Furthermore, if $u \in A$ and $v \notin A$, then $(A \setminus \{u\}) \cup \{v\}$ also resolves $G$.
\end{lem}

\begin{thm} \cite{BH} \label{t6.1}
	If $A, B, C$ and $D$ are square matrices of the
	same order and $D$ is invertible. Then \[\begin{vmatrix}
	A & B  \\
	C & D \end{vmatrix} = \begin{vmatrix}D 	\end{vmatrix}. \begin{vmatrix} A - BD^{-1}C \end{vmatrix}. \]
\end{thm}

\begin{thm} \cite{Cvetkovi{\'c}(2010)} \label{t6.2}
	For any vertex $u$ of a connected graph $G$, we have
	$$ \lambda_1(G-u) <  \lambda_1(G).$$
\end{thm}

\begin{thm} \cite{Cvetkovi{\'c}(2010)} \label{t6.3}
	Let $A, B$ be two $n \times n$ hermitian matrices. Then
	
	$$ \lambda_{r}(A+B) \leq  \lambda_{s}(A) + \lambda_{r-s+1}(B), \text{ where } n \geq r \geq s \geq 1. $$
\end{thm}

\begin{lem}
	The metric dimension of $P(G(n))$ is $\psi(P(G(n))=  2^{n} -3.$
	
\end{lem}

\noindent{\textbf{Proof.}} Note that, $P(n) \setminus \{e\}$ and $H(n)$ are twin sets in $P(G(n))$ of order $2^{n-1}-1$ and $2^{n-1}$ respectively.
By using the argument of Lemma \ref{l6.2}, $P(n) \setminus \{e\}$ and $H(n)$ contains at least $(2^{n-1}-1)-1$ and $2^{n-1}-1$ number of elements. Thus the metric dimension of $P(G(n))$ is at least  $2^{n}-3$, i.e., $\psi(P(G(n)) \geq 2^n - 3$. Also,  $B = \{P(n) \setminus \{0,1\}\} \cup \{H(n) \setminus {2^{n-1}}\}$ is the resolving set for $P(G(n))$ of order $2^n-3$, thus $\psi(P(G(n)) \leq 2^n - 3$. Hence, $\psi(P(G(n)) = 2^n - 3$. \hfill $\Box$

\begin{lem}
	Let $P(G(n))$ be the power graph of $G(n)$. Then
	$\psi(P(G(n)), x) = x^{2^n} + 2^nx^{2^n-1} + (2^{2n-2} + 2^{n-1} -1)x^{2^n-2} + (2^{2n-2}-2^{n-1})x^{2^n-3}.$
\end{lem}

\noindent{\textbf{Proof.}} Since the metric dimension of $P(G(n))$ is $2^n-3,$ we need to find only $r_{2^n-3},r_{2^n-2},r_{2^n-1}, \ \text{and} \ \\ r_{2^n}.$ As each resolving set can omit at most one element from each twin sets $P(n) \setminus \{e\}$ and $H(n)$, so $r_{2^n-3}$ must contains $2^{n-1}-2$ elements from $P(n) \setminus \{e\}$ and $2^{n-1}-1$ elements from  $H(n).$ Thus all the resolving sets for $P(G(n))$ of order $2^n-3$ are
$$r_{2^n-3} = {2^{n-1} \choose 2^{n-1}-1} {2^{n-1}-1 \choose 2^{n-1}-2} = 2^{n-1}(2^{n-1}-1).$$

Let $A$ be a resolving set of $P(G(n))$ of order $2^n-2$ and $u,v \in V(P(G(n))$ such that $u,v \notin A$. Since each resolving set can omit at most one element from each twin sets $P(n) \setminus \{e\}$ and $H(n)$, we have either $u=e$ and $v \in H(n)$, either $u=e$ and $v \in P(n) \setminus \{e\}$, either $u \in H(n)$ and $v \in P(n) \setminus \{e\}$. Thus, all the resolving sets for $P(G(n))$ of order $2^n-2$ are

$$r_{2^n-2} = {2^{n-1} \choose 2^{n-1}-1} {2^{n-1}-1 \choose 2^{n-1}-1} \ \ \text{or} \ \ {2^{n-1} \choose 2^{n-1}} {2^{n-1}-1 \choose 2^{n-1}-2} \ \ \text{or} \ \ {2^{n-1} \choose 2^{n-1}-1} {2^{n-1}-1 \choose 2^{n-1}
	-2}.$$

Hence, $r_{2^n-2} = 2^{2n-2} + 2^{n-1} -1.$ Also, the resolving set having $2^n-1$ order can be selected into $2^n$ possibly distinct ways. Thus, $r_{2^n-1} = 2^n$. Clearly, $r_{2^n} = 1$. This gives,  

$$ \psi(P(G(n))), x) = x^{2^n} + 2^nx^{2^n-1} + (2^{2n-2} + 2^{n-1} -1)x^{2^n-2} + (2^{2n-2}-2^{n-1})x^{2^n-3}.  \hspace{2cm }\hfill \Box $$

\begin{lem}
	The characteristic polynomial of $P(G(n))$ is $$x^{2^{n-1}-1}(1+x)^{(2^{n-1}-2)}[x^3 + x^2(2-2^{n-1})x^2 - (1-2^n)x + 2^{2n-2}-2^n].$$
\end{lem}

\noindent{\textbf{Proof.}}  The adjacency matrix $A_{G(n)}$ of $P(G(n))$ is $
\begin{bmatrix}

B_{2^{n-1}} & C_{2^{n-1}}  \\
C^T_{2^{n-1}} & O_{2^{n-1}} 
\end{bmatrix}$, where

\[B_{2^{n-1}} =
\begin{bmatrix}

0 & 1 & 1 & \cdots & 1 \\
1 & 0 &  1 &\cdots& 1\\
1 & 1 & 0 &\cdots& 1\\
\vdots & \vdots & \vdots &\ddots & \vdots\\
1 & 1 & 1 & \cdots & 0
\end{bmatrix} \text{and} \ C_{2^{n-1}} =
\begin{bmatrix}

1 & 1 & 1 & \cdots & 1 \\
0 & 0 &  0 &\cdots& 0\\
0 & 0 & 0 &\cdots& 0\\
\vdots & \vdots & \vdots &\ddots & \vdots\\
0 & 0 & 0 & \cdots & 0
\end{bmatrix}
\]

So, by Theorem \ref{t6.1}, the characteristic polynomial $p(x)$ of $P(G(n))$ is

\[p(x) =
\begin{vmatrix}
xI_{2^{n-1}}-
B_{2^{n-1}} & -C_{2^{n-1}}  \\ 
-C^T_{2^{n-1}} & xI_{2^{n-1}}
\end{vmatrix} = \begin{vmatrix}
xI_{2^{n-1}}
\end{vmatrix}.\begin{vmatrix}
xI_{2^{n-1}}-
B_{2^{n-1}}-C_{2^{n-1}}\frac{I_{2^{n-1}}}{x}C^T_{2^{n-1}}
\end{vmatrix}  
\]
Now,
\begin{align*}
\begin{vmatrix}
xI_{2^{n-1}}-
B_{2^{n-1}}-C_{2^{n-1}}\frac{I_{2^{n-1}}}{x}C^T_{2^{n-1}}
\end{vmatrix}
&  = \begin{vmatrix}
x - \frac{2^{n-1}}{x} & -1 & -1 & \cdots & -1 \\
-1 & x &  -1 &\cdots& -1\\
-1 & -1 & x &\cdots& -1\\
\vdots & \vdots & \vdots &\ddots & \vdots\\
-1 & -1 & -1 & \cdots & x
\end{vmatrix}_{2^{n-1} \times 2^{n-1}} \hspace{-.5cm} (\text{$R_1 \rightarrow R_1-R_2$})  \\
& = \begin{vmatrix}
x - \frac{2^{n-1}}{x}+1 & -1-x & 0 & \cdots & 0 \\
-1 & x &  -1 &\cdots& -1\\
-1 & -1 & x &\cdots& -1\\
\vdots & \vdots & \vdots &\ddots & \vdots\\
-1 & -1 & -1 & \cdots & x
\end{vmatrix}_{2^{n-1}\times 2^{n-1},} \\
& = (x - \frac{2^{n-1}}{x}+1).\begin{vmatrix}
x &  -1 &\cdots& -1\\
-1 & x &\cdots& -1\\
\vdots & \vdots &\ddots & \vdots\\
-1 & -1 & \cdots & x
\end{vmatrix}_{2^{n-1}-1\times 2^{n-1}-1,} \\
& + (1+x).\begin{vmatrix}
-1 &  -1 &\cdots& -1\\
-1 & x &\cdots& -1\\
\vdots & \vdots &\ddots & \vdots\\
-1 & -1 & \cdots & x
\end{vmatrix}_{2^{n-1}-1\times 2^{n-1}-1,} \\ & = (x - \frac{2^{n-1}}{x}+1)D + (1+x)E,
\end{align*} where

$D = \begin{vmatrix}
x &  -1 &\cdots& -1\\
-1 & x &\cdots& -1\\
\vdots & \vdots &\ddots & \vdots\\
-1 & -1 & \cdots & x
\end{vmatrix}_{2^{n-1}-1\times 2^{n-1}-1,} \text{ and } E = \begin{vmatrix}
-1 &  -1 &\cdots& -1\\
-1 & x &\cdots& -1\\
\vdots & \vdots &\ddots & \vdots\\
-1 & -1 & \cdots & x
\end{vmatrix}_{2^{n-1}-1\times 2^{n-1}-1,}$

\bigskip
Apply $C_1 \rightarrow C_1 + C_2 + \cdots + C_{2^{n-1}-1}$ in $D$, we have

\[ \begin{vmatrix}
x -(2^{n-1}-2)&  -1 &\cdots& -1\\
x -(2^{n-1}-2) & x &\cdots& -1\\
\vdots & \vdots &\ddots & \vdots\\
x -(2^{n-1}-2) & -1 & \cdots & x
\end{vmatrix}_{2^{n-1}-1\times 2^{n-1}-1,} = (x -(2^{n-1}-2)).\begin{vmatrix}
1 &  -1 &\cdots& -1\\
1 & x &\cdots& -1\\
\vdots & \vdots &\ddots & \vdots\\
1 & -1 & \cdots & x
\end{vmatrix}_{2^{n-1}-1\times 2^{n-1}-1,} \]

Apply $R_2 \rightarrow R_2 - R_1$, $R_3 \rightarrow  R_3 - R_1, \ldots, R_n \rightarrow  R_n - R_1$

\[ (x -(2^{n-1}-2)) \begin{vmatrix}

1 &  -1 &\cdots& -1\\
0 & x+1 &\cdots& 0\\
\vdots & \vdots &\ddots & \vdots\\
0 & 0 & \cdots & x+1
\end{vmatrix}_{2^{n-1}-1\times 2^{n-1}-1,} = (x -(2^{n-1}-2))(x +1)^{(2^{n-1}-2)}) \]

Now, apply $R_2 \rightarrow  R_2 - R_1$, $R_3 \rightarrow  R_3 - R_1, \ldots, R_n \rightarrow  R_n - R_1$ in $E$, we have

\[ \begin{vmatrix}
-1 &  -1 &\cdots& -1\\
0 & 1+x &\cdots& -1\\
\vdots & \vdots &\ddots & \vdots\\
0 & 0 & \cdots & 1+x
\end{vmatrix}_{2^{n-1}-1\times 2^{n-1}-1,} \\
= -(1+x)^{2^{n-1}-2}.  \] Thus,

\begin{align*}
\begin{vmatrix}
xI_{2^{n-1}}-
B_{2^{n-1}} & -C_{2^{n-1}}  \\ 
C^T_{2^{n-1}} & xI_{2^{n-1}}
\end{vmatrix} 
& = x^{2^{n-1}}\bigg[(
x - \frac{2^{n-1}}{x}+1)(x -(2^{n-1}-2))((x +1)^{(2^{n-1}-2)}) \\ 
& -  (1+x)(1+x)^{(2^{n-1}-2)}\bigg]\\
& = x^{2^{n-1}-1}(1+x)^{(2^{n-1}-2)}\bigg[x^3 + x^2(2-2^{n-1})x^2- (1-2^n)x + 2^{2n-2}-2^n\bigg]. 
\end{align*} \hfill $\Box$

\begin{lem}
The spectral radius of  $P(G(n))$ is
$$ \lambda_1(P(Z_{2^{n-1}})) < \lambda_1(P(G(n)) \leq \lambda_1(P(Z_{2^{n-1}}) + \sqrt{{2^{n-1}}}$$
\end{lem}

\noindent{\textbf{Proof.}}  The adjacency matrix of  $P(G(n))$ is $A = D +E$, where 
$D =
\begin{bmatrix}

B_{2^{n-1}} & O_{2^{n-1}}  \\
O_{2^{n-1}} & O_{2^{n-1}} 
\end{bmatrix}$, and $E =
\begin{bmatrix}

O_{2^{n-1}} & C_{2^{n-1}}  \\
C^T_{2^{n-1}} & O_{2^{n-1}} 
\end{bmatrix}$ are hermitian matrices. Since $V(P(G(n))) = P(n) \cup H(n)$, where $P(n)$ is a cyclic group $Z_{2^{n-1}}$, the adjacency matrix of power graph of $P(n)$ is $B_{2^{n-1}}$. This gives,  $\lambda_1(D) = \lambda_1(P(Z_{2^{n-1}})) $. By Theorem \ref{t6.3}, we have 

$$ \lambda_1(P(G(n))) \leq \lambda_1(D) + \lambda_1(E) = \lambda_1(P(Z_{2^{n-1}})) + \lambda_1(E) .$$

Now,

$$\begin{vmatrix}
xI_{2^{n-1}}-E
\end{vmatrix} =
\begin{vmatrix}

xI_{2^{n-1}} & -C_{2^{n-1}}  \\
-C^T_{2^{n-1}} & xI_{2^{n-1}} 
\end{vmatrix}$$

First apply $R_1 \rightarrow xR_1$, and then $R_1 \rightarrow R_1+ R_{2^{n-1}+1}$, $R_1 \rightarrow R_1 + R_{2^{n-1}+2}, \ldots R_1 \rightarrow R_1 + R_{2^{n}}$, we get

$$\begin{vmatrix}
xI_{2^{n-1}}-E
\end{vmatrix} = \begin{vmatrix}

(x^2-2^{n-1})I_{2^{n-1}} & O_{2^{n-1}}  \\
-C^T_{2^{n-1}} & xI_{2^{n-1}}  
\end{vmatrix} = 
(x^2-2^{n-1})^{2^{n-1}}.x^{2^{n-1}}$$

Thus, the eigen values of $E$ are $  \underbrace{\pm \sqrt{{2^{n-1}}},\ldots,\pm \sqrt{{2^{n-1}}}}_{2^{n-1}\ \text{times}},\underbrace{0,\ldots,0}_{2^{n-1}\ \text{times}}.$  Therefore $\lambda_1(E) = \sqrt{{2^{n-1}}}.$

\smallskip
Also, by Theorem \ref{t6.3}, $\lambda_1(P(Z_{2^{n-1}})) < \lambda_1(P(G(n)).$ Combining all the inequalities, we have

$$ \lambda_1(P(Z_{2^{n-1}})) < \lambda_1(P(G(n)) \leq \lambda_1(P(Z_{2^{n-1}}) + \sqrt{{2^{n-1}}}.  $$  $ \hfill \Box$

\section{Distant and detour distant properties}

In this section, we compute closure, interior, distance degree sequence of the $P(G(n))$. 
\begin{lem}
	Let $P(G(n))$ be the power graph on $G(n)$. Then,
	
	\begin{enumerate}
		\item 
		\begin{equation*}
		ec_D(u)= 
		\begin{cases}
		2^{n-1}; & \text{$ \text{if} \ u \in P(n) \setminus \{e\},$}\\
		2^{n-1}-1; & \text{$ \text{if} \ u = e,$}\\	
		2^{n-1}; & \text{$ \text{if} \ u \in H(n),$}
		
		\end{cases}       
		\end{equation*}
		
		\item $rad_D(P(G(n))) = 2^{n-1} -1,$
		
		\item $dia_D(P(G(n))) = 2^{n-1}.$
		
	\end{enumerate}
\end{lem}

\noindent{\textbf{Proof:}} Clearly, every two distinct vertices of $P(n)$ are adjacent, but every two distinct vertices of $H(n)$ are non-adjacent. Also, $e$ is adjacent to all the vertices of $P(G(n))$ and $u$ is non-adjacent with $v$ for every $u \in P(n) \setminus \{e\}$ and $v \in H(n).$ If $u \in P(n) \setminus \{e\}$, then there is a $u$-$v$ path having detotar length $2^{n-1}-1$ for every $v \in P(n)$ and $u$-$w$ path having detotar length $2^{n-1}$ for every $w \in H(n)$. Thus, $ec_D(u) = 2^{n-1}$ for $u \in k_n \setminus \{e\}$. If $u =e$, then there is a $u$-$x$ path having detotar length $2^{n-1} - 1$ for every $v \in P(n) \setminus \{e\}$ and a $u$-$y$ path having detotar length $1$ for every $v \in H(n)$. Thus, $ec_D(u) = 2^{n-1} - 1$ for $u = e$. Now, if $u \in H(n)$, then there is a $u$-$z$ path having detotar length $2^{n-1}$ for every $v \in P(n) \ \setminus \{e\}$ and $u$-$e$ path having detotar length $1$. Thus, $ec_D(u) = 2^{n-1}$ for $u \in H(n)$. Since, the smallest and largest detotar eccentricity of $P(G(n))$ is $2^{n-1}-1$ and $2^{n-1}$ respectively, it follows that $rad_D(P(G(n))) = 2^{n-1} -1$ and $dia_D(P(G(n))) = 2^{n-1}.$ \hfill $\Box$

\begin{lem}\label{l5.1}
	Let $P(G(n))$ be the power graph on $G(n)$. Then, the distance degree sequence and detour distance degree sequence of the power graph $P(G(n))$, is given
	by
	
	\begin{enumerate}
		\item 
		$dds(P(G(n))) = ((1,2^n-1),(1,2^{n-1}-1,2^{n-1})^{2^{n-1}-1},(1,1,2^n-2)^{2^{n-1}}),$ 
		
		\item $dds_D(P(G(n))) = ((1,2^{n-1},\underbrace{0,\ldots,0}_{2^{n-1}-3\ \text{times}},2^{n-1}-1),(1,\underbrace{0,\ldots,0}_{2^{n-1}-2\ \text{times}},2^{n-1}-1,2^{n-1})^{2^{n-1}-1}, \linebreak (1,1,2^{n-1}-1,\underbrace{0,\ldots,0}_{2^{n-1}-3\ \text{times}},2^{n-1}-1)^{2^{n-1}}$.
		
	\end{enumerate}
	
	\noindent{\textbf{Proof:}} For $u \in V(P(G(n))) \setminus \{e\}$, we get $ec(u) = 2$ and by Lemma \ref{l5.1}, $ec_D(u) = 2^{n-1}$. Next, for $u = e$, we get $ec(u) = 1$ and by Lemma \ref{l5.1}, $ec_D(u) = 2^{n-1} - 1$. Therefore
	
	\begin{equation*}
	dd(u)= 
	\begin{cases}
	(1, 2^n-1); & \text{$ \text{if} \ u = e,$}\\
	(1, 2^{n-1} -1, 2^{n-1}); & \text{$ \text{if} \ u = P(n) \setminus \{e\},$}\\	
	(1,1,2^n-2); & \text{$ \text{if} \ u \in H(n),$}
	
	\end{cases}       
	\end{equation*}

	\begin{equation*}
	dds_D(u)= 
	\begin{cases}
	(1, 2^{n-1},\underbrace{0,\ldots,0}_{2^{n-1}-3\ \text{times}},2^{n-1}-1); & \text{$ \text{if} \ u = e,$}\\
	(1,\underbrace{0,\ldots,0}_{2^{n-1}-2\ \text{times}},2^{n-1}-1,2^{n-1}); & \text{$ \text{if} \ u = k_n \setminus \{e\},$}\\	
	(1,1,2^{n-1}-1,\underbrace{0,\ldots,0}_{2^{n-1}-3\ \text{times}},2^{n-1}-1); & \text{$ \text{if} \ u \notin k_n,$}
	
	\end{cases}       
	\end{equation*}
	
	Since,  $|\{e\}|$,  $|P(n)\setminus \{e\}|, \text{and} \ |H(n)|$ are  $1, 2^{n-1}-1, \text{and} \ 2^{n-1}$ respectively. Therefore, 
	
	$$dds(P(G(n))) = ((1,2^n-1),(1,2^{n-1}-1,2^{n-1})^{2^{n-1}-1},(1,1,2^n-2)^{2^{n-1}}),$$ 
	
	$dds_D(P(G(n))) = ((1,2^{n-1},\underbrace{0,\ldots,0}_{2^{n-1}-3\ \text{times}},2^{n-1}-1),(1,\underbrace{0,\ldots,0}_{2^{n-1}-2\ \text{times}},2^{n-1}-1,2^{n-1})^{2^{n-1}-1}, \linebreak  (1,1,2^{n-1}-1,\underbrace{0,\ldots,0}_{2^{n-1}-3\ \text{times}},2^{n-1}-1)^{2^{n-1}}$.   \hfill $\Box$
	
\end{lem}
\smallskip
Now, we compute interior and closure of $P(G(n))$ by using below results.

\begin{lem}\label{l16} (\cite{cha(2006)})
	Let $G$ be a connected graph and $u \in V(G)$. Then $u$ is said to be a boundary vertex of each vertex other than $u$ if and only if $u$ is a complete vertex of $G$.
\end{lem}

\begin{lem}\label{l17} (\cite{cha(2006)})
	Let $G$ be a connected graph and $u \in V(G)$. Then $u$ is said to be a boundary vertex of $G$ if and only if $u$ is not an interior
	vertex of $G$
\end{lem}

\begin{lem}
	Let $P(G(n))$ be the power graph of $G(n)$. Then
	$$Int(P(G(n))) = Ce(P(G(n))) \ \text{and} \ Cl(P(G(n))) = P(G(n)).$$
\end{lem}

\noindent{\textbf{Proof.}}

\begin{equation*}
N(u)= 
\begin{cases}
P(n) \setminus \{u\}; & \text{$ \text{if} \ u \in P(n) \setminus \{e\},$}\\
G(n) \setminus \{e\}; & \text{$ \text{if} \ u = e,$}\\	
e; & \text{$ \text{if} \ u \in H(n),$}

\end{cases}       
\end{equation*}

Note that, the subgraphs induced by $\{e\}$ and $P(n) \setminus \{e\}$ are complete and the induced subgraph of $P(G(n)) \setminus \{e\}$ is disconnected. This gives that, every vertex $u$ of $P(G(n)) \setminus \{e\}$ is a complete vertex. Hence, by Lemma \ref{l16}, $u$ is a boundary vertex. So, by Lemma \ref{l17}, $e$ is an interior vertex of the graph $P(G(n))$. Clearly $Ce(P(G(n)) = \{e\}$. Hence, we conclude that $Int(P(G(n)) = Ce(P(G(n))$.

\smallskip

Let $u, v \in P(G(n))$. If $u \in P(n) \setminus \{e\}$ and $v \in H(n)$, then $deg(u) = 2^{n-1}-1$ and $deg(v) = 1$. So, $deg(u) + deg(v) = 2^{n-1}-1 +1 = 2^{n-1} < 2^n.$ Also, if  $u \in H(n)$ and $v \in H(n)$, then $deg(u) = 1$ and $deg(v) = 1$. This gives that $deg(u) + deg(v) = 1 +1 < 2^n.$ Hence, $Cl(P(G(n))) = P(G(n)).$ \hfill$\Box$

\section{Discussion}


Motivated by the earlier study of power graph for a group, here, we have described the notion of such graph for a gyro-group for the first time. Using this, we have obtained some of its combinatorial properties for the power graph of gyrogroup $G(n)$, defined in Section \ref{s3}. From the study, we have seen that finite gyro-commutative gyrogroups with isomorphic power graphs need not be isomorphic. This raised the following question: 

\begin{ques}
	If $G_1, G_2$ are two non-degenrate gyrogroups, then under what condition $P(G_1) \cong P(G_2) \Longrightarrow G_1 \cong G_2$?
\end{ques}

The properties of the power graph of gyrogroups, in comparison to the properties of the power graph of groups, have not yet been thoroughly investigated. In this context, we would like to finalise this discussion by posing some questions.

Cameron and Ghosh \cite{CG} shown that the only finite group whose automorphism group is the same as its power graph is the Klein group. Thus it is natura to ask:

\begin{ques}
	Is there any non-degenerate gyrogroup $G$ whose automorphism group is the same as its power graph?
\end{ques}

Ivy Chakrabarty e.t. al \cite{Chakrabarty(2009)} shown that for a finite group $G$, $P(G)$ is complete if and only if $G$ is a cyclic group of order 1 or $p^m$, for some prime $p$ and for some natural number $m$. This leads to the following question:

\begin{ques}
	In which condition power graph $P(G)$ of a non-degenerate gyrogroup $G$ is complete?
\end{ques}

Mirzargar et al. \cite{MAA} shown that the power graph $P(Z_{p^n})$ has the maximum number of edges among all power graphs of $p$-groups of order $p^n$. Also, they conjectured that the power graph $P(Z_n)$ has the maximum number of edges among all power graphs of finite groups of order $n$.  This arises the following question:

\begin{ques}
	For which non-degenerate gyrogroup $G$, $P(G)$ has the maximum number of edges among all power graphs of finite non-degenerate gyrogroup of order $n$?
\end{ques}

\section{Acknowledgment}

The authors thank Professor Teerapong Suksumran for his valuable suggestions, discussions and constant support. The first author is thankful to the Ministry of Human Resource Development (MHRD) New Delhi,
India for financial support.

\newpage
\end{document}